\documentclass[12pt,fleqn,emtex]{article}

\usepackage{graphicx}
\usepackage{framed}
\usepackage{times}
\usepackage{mathptmx}

\usepackage[german, american]{babel}  %Deutsch
\usepackage[utf8]{inputenc}

\usepackage{amsfonts}       %
\usepackage{amsmath} %[leqno]=Gleichungen werden links numeriert;
\usepackage{amssymb}        %
\usepackage{amstext}        %
\usepackage{amsthm}         %theorem-umgebungen
\usepackage[matrix,curve]{xypic}

\usepackage{textcomp}
\usepackage{epsfig}
\psfigdriver{dvips}
\usepackage{latexsym}
\usepackage[matrix,curve]{xypic}
\usepackage{rotating}
\rotdriver{dvips}

\newcommand{\Z}{\mathbb{Z}}
\newcommand{\Q}{\mathbb{Q}}

\newcommand{\C}{\mathbb{C}}

\newtheorem{Lemma}{Lemma}
\newtheorem{Theorem}{Theorem}
\newtheorem{Corollary}{Corollary}
\newtheorem{Proposition}{Proposition}

\DeclareMathSizes{12}{11}{8}{5} 
\usepackage{xypic,dsfont}

\begin{document}

\title{On the Rigidity of BN-sheaves}
\author{R. Weissauer}

\maketitle
\thispagestyle{empty}

%\section*{Introduction}
\bigskip\noindent

Let $X$ be an abelian variety over
an algebraically closed field $k$, where we assume that either $k=\C$ or $k$ is the algebraic closure of a finite field.  In [BN] and [W2] we considered the convolution
product $K*L$ for complexes $K$ and $L$ in the bounded 
derived category $D_c^b(X,\Lambda)$, where the coefficient field $\Lambda$
is either $\C$
for $k=\C$ or $\Lambda=\overline \Q_l$. The convolution product 
is defined by the group law $a: X\times X\to X$ of the abelian 
variety $X$, as the derived direct image complex $K*L=Ra_*(K\boxtimes L)$. 
This convolution product makes $(D_c^b(X,\Lambda),*)$ into a rigid triangulated symmetric monoidal $\Lambda$-linear category; its unit object is the skyscraper sheaf $\delta_0$ concentrated at zero. 
For further details we refer
to [W2]. For a complex $K$ let $D(K)$ denote its Verdier dual and $K^\vee = (-id_X)^*D(K)$ its rigid dual.

\medskip
For our considerations, the decomposition theorem and the hard Lefschetz theorem for perverse sheaves are essential perequisites. For this we specify
a full  $\Lambda$-linear suspended tensor subcategory $({\bf D},*) \subseteq (D_c^b(X,\Lambda),*)$ 
as in [KrW, example 6], so that among others objects in ${\bf D}$ are semisimple, the decomposition theorem holds and also the hard Lefschetz theorem. In particular the perverse cohomology functors ${}^pH^i(K) \in {\bf P}$ are defined for $K$ in $\bf D$ where ${\bf P} \subset {\bf D}$ is an abelian subcategory
of perverse sheaves defined by a perverse $t$-structure on ${\bf D}$ with core ${\bf P}$.  If we speak of perverse sheaves on $X$,  we always mean objects in this category ${\bf P}$. Let $e$ denote the projector $e:{\bf D}\to {\bf P}$, then $e[n]$ is the projector to ${\bf P}[n]$.
The categories $\bf D$ and ${\bf P}$ are stable under twists $K \mapsto K_\chi = K\otimes_{\Lambda_X} L_\chi$ with respect to local systems $L_\chi$ defined by the characters $\chi$ of the fundamental group $\pi_1(X,0)$ of $X$ in the sense of [KrW]
and stable under $K \mapsto T^*_x(K)$ for translations $T_x(y)=y+x$ with respect to closed points $x\in X$. 

\medskip
{\bf Evaluation morphisms}. 
We now discuss properties of the suspended symmetric monoidal rigid  
$\Lambda$-linear  tensor category $({\bf D},*)$ with the tensor product $*$ defined by the convolution (see [BN] and [KrW], also for the notations used). ${\bf D}$ is a Krull-Schmidt category, i.e. an additive category for which every object decomposes into a finite direct sum of objects with local endomorphism rings. A Krull-Schmidt category is idem-complete, and its objects are indecomposable if and only if their endomorphism ring is local. Any isomorphism $\bigoplus_{i=1}^n K_i \cong \bigoplus_{j=1}^m L_j$ for indecomposable objects $K_i,L_j$ implies $n=m$ and $L_j \cong K_{\sigma(i)}$ for a permutation $\sigma$. In fact, any object in $\bf D$ decomposes into a finite direct sum of objects not necessarily equal, but isomorphic an object of the form $P[n]$ for irreducible $P\in \bf P$ and $n\in \Z$. Since $End_{\bf P}(P[n]) \cong \Lambda \cdot id_{P[n]}$, the category is Krull-Schmidt, and
 the indecomposable
objects $K$ in $\bf D$ have the form $K=P[n]$ for irreducible $P\in \bf P$ and $n\in \Z$. 
By rigidity [W2], for any $K$ in ${\bf D}$ there exists a coevaluation morphism
$$  coev_K: \delta_0 \to K*K^\vee \ ,$$  
corresponding to the identity $id_K$ via $Hom_{D_c^b(X,\Lambda)}(K,K)= Hom_{D_c^b(X,\Lambda)}(\delta_0,K*K^\vee)$. Similarly one has the 
evaluation morphism $$  eval_K:  K^\vee *K \to \delta_0 \ ,$$
so that the composition of the induced morphisms $(id_K*eval_K)\circ ass \circ (coev_K*id_K)$
%resp. $(eval_K*id_{K^\vee})\circ ass \circ (id_{K^\vee} *coev_K)$ 
$$ K = \delta_0 * K \to (K*K^\vee)*K \to K*(K^\vee*K)\to K*\delta_0 = K$$
%$$ K^\vee = K^\vee *\delta_0 \to K^\vee *(K*K^\vee) \to (K^\vee *K)*K^\vee\to \delta_0 *K^\vee = K^\vee$$
is the identity morphism $id_K:K\to K$. There is a similar dual identity for $K^\vee$. 
%respectively the identity morphism of $K^\vee$.   

\medskip
{\bf Remark}.   $(K^\vee,eval_K)$ attached to an object $K$ is unique up to isomorphism (see [CT,p.120]). We use this together with the following simple facts (see [D,1.15]).

\medskip
a) Suppose $K$ is a retract of $L$ defined by an idempotent $e\in End_{\bf D}(L)$ admitting a direct sum (i.e. biproduct) decomposition. Then for $K^\vee$, considered as a retract of $L^\vee$ defined via the dual idempotent $e^\vee$, this gives a retract $\iota: K^\vee * K \hookrightarrow L^\vee * L$ so that $eval _K =eval_L \circ \iota$ holds.

\medskip 
b) For $K=A\oplus B$ the evaluation $eval_K$ is obtained as $eval_A + eval _B$,
using the projection $K^\vee*K \twoheadrightarrow (A^\vee*A) \oplus (B^\vee* B)$.

\medskip
c) For $K=A*B$, using $K^\vee*K = (A*B)^\vee *(A*B)\cong (A^\vee*A)*(B^\vee*B)$, the evaluation morphism of $K^\vee*K\to \delta_0$ is obtained as the tensor product
$eval_K = eval_A * eval_B$.  

\medskip
The symmetry constraints of the tensor category define
isomorphisms $$S: K*K^\vee \cong K^\vee *K$$ such that the composed
morphism $eval_K \circ S \circ coev_K$
$$ \delta_0 \to K*K^\vee \to K^\vee * K \to \delta_0 \ ,$$
considered as an element of $Hom_{\bf D}(\delta_0,\delta_0)=\Lambda$,
is the multiplication with the categorial dimension of $K$; in our case the categorial
dimension is the Euler characteristic 
$\chi(K) = \sum_i (-1)^i \dim_\Lambda(H^i(X,K))$
of the complex $K$.

\medskip
For a morphism $\rho:K\to L$ the transposed
morphism $\rho^\vee: L^\vee \to K^\vee$ is defined as
$(eval_L*id_{K^\vee})\circ(id_{L^\vee}*\rho* id_{K^\vee})\circ(
id_{L^\vee}*coev_{K})$. Together with $K \mapsto K^\vee$ this induces a tensor equivalence
with the opposite category so that $(K^\vee)^\vee\cong K$ and $(\rho^\vee)^\vee=\rho$
in the sense of [CT,2.5]. There exists an isomorphism $\varphi: (K^{\vee\vee}*K^\vee)^\vee \cong K*K^\vee$ such that $coev_{K} = \varphi \circ (eval_{K^\vee})^\vee  $. Using the duality isomorphisms $d_{K,L}:  L^\vee *K^\vee \to (K*L)^\vee$
defined by $(eval_L*id_{(K*L)^*})\circ(id_{L^\vee}*eval_K*id_L*id_{(K*L)^\vee}) \circ
(id_{L^\vee}*id_{K^\vee}*coev_{K*L})$, more concretely one can show $(coev_{K})^\vee \circ D= eval_{K^\vee}$ 
for $D=d_{K,K^\vee}$.

\medskip
{\bf Monoidal components}.  By the decomposition theorem 
$K^\vee*K $ is semisimple for $K\in {\bf D}$. Hence $K^\vee*K  = \bigoplus_i \ {}^p H^{i}(K^\vee*K)[-i]$,
and any  ${}^p H^{i}(K^\vee*K)$ decomposes into a direct sum $\bigoplus P_\nu^{i}$
of irreducible perverse sheaves $P_\nu^{i}$. 
Using this decomposition, the evaluation can be written as 
a sum $ eval_K = \sum_{\nu,i} ev_{\nu,i} $ with morphisms $ev_{\nu,i} \in Hom_{\bf D}(P_{\nu}^{i}[-i], \delta_0)$.  Since for irreducible $K$ $$Hom_{\bf D}(K,K)= Hom_{\bf D}(K^\vee*K,\delta_0)$$
has dimension one, there exists
a unique exponent $i=\nu_K$ and a unique simple perverse constituent ${\cal P}_K$ of 
${}^p H^{i}(K*K^\vee)$ such that $eval_K$ factorizes over ${\cal P}_K[-i]$. All the other morphisms $ev_{\nu,i}$ are zero. This gives a commutative diagram
 $$ \xymatrix{ \ \ K^\vee*K \ \ \ar@<1ex>@{->>}[dd]^p \ar[rd]^{eval_K} &  \cr
& \ \ \delta_0 \cr
\ \ {\cal P}_K[-\nu_K]\ \ \ar@<1ex>@{^{(}->}[uu]^\iota \ar[ru]_\epsilon & } \ $$
where $p\circ \iota = id$ is the identity morphism. In the following,
arrows $\hookrightarrow$ and $\twoheadrightarrow$ always split monomorphisms $\iota $ and the corresponding projections $p$ obtained from direct sum decompositions, which makes sense in our  $\Lambda$-linear tensor category $({\bf D},*)$.  However, for convenience, we reserve these symbols for retracts associated to idemponents $\iota\circ p$ that commute with all idemponents $e[n]$. Put briefly, this assures that the correspond to decompositions into direct sums of translates of perverse sheaves. This property is preserved by functors $Rf_*$, hence by the convolution product. For an arbitrary  rigid symmetric monoidal $\Lambda$-linear tensor category $eval_K=0$ implies $id_K=0$ and hence $K=0$,
this shows $\epsilon\neq 0$. 

\medskip
For an irreducible perverse sheaf $K$ the distinguished irreducible component ${\cal P}_K$  will be called the {\it monoidal component} of the irreducible perverse sheaf $K$, and $\nu_K$ its {\it degree}, and in the case  $\nu_K>0$ the perverse sheaf ${\cal P}_K$ will be called a monoidal perverse sheaf or {\it monoid} on $X$.  Concerning this, notice that the degree always satisfies $\nu_K\geq 0$. 

\medskip
This follows from the perverse vanishing conditions $$ Hom_{\bf D}(M,N[r])=0$$ for $M,N\in {\bf P}$ and $r<0$, applied for the objects ${\cal P}_K$  and
$\delta_0$ in ${\bf P}$. 

\medskip
From the definition of ${\cal P}_K[\nu_K]$ and the existence of the symmetry isomorphism
$S: K*K^\vee \cong K^\vee*K$, it is clear that ${\cal P}_{K^\vee} \cong {\cal P}_K$ and $\nu_{K^\vee}=\nu_K$. Therefore passing to the dual, using 
$coev_{K} = \varphi \circ (eval_{K^\vee})^\vee  $ for some isomorphism $\varphi: (K^{\vee\vee}*K^\vee)^\vee \cong K*K^\vee$ and ${\cal P}_K[-\nu_K]^\vee \cong {\cal P}_K^\vee [+\nu_K]$,
we obtain a commutative diagram
$$ \xymatrix{  & \ \ K*K^\vee \ \ \ar@<1ex>@{->>}[dd]^{p'}   \cr
\delta_0 \ar[ru]^{coev_K} \ar[rd]_\sigma \ \ &  \cr
& \ \ {\cal P}_K^\vee[+\nu_K]\ \ \ar@<1ex>@{^{(}->}[uu]^{\iota'}   } \ $$

Using the perverse vanishing condition for morphisms and the adjunction formulas \label{adformeln} $$Hom_{\bf D}(K,\delta_0[n]) = {\cal H}^{-n}(K)_0^*$$ and $Hom_{\bf D}(K*L,\delta_0)= Hom_{\bf D}(L,K^\vee)$ for $K,L\in {\bf D}$ 
it follows that ${\cal H}^{>0}(K*L)=0$ holds for perverse sheaves $K,L \in {\bf P}$.
Hence the following assertions 5, 6, 7 and 9 of lemma \ref{1} are an immediate consequence,
in view of the hard Lefschetz theorem. 

\begin{Lemma} \label{1} Suppose $K \in {\bf P}$ is an irreducible perverse sheaf on $X$, then
\begin{enumerate}
\item ${\cal P}_{K}$ is irreducible and ${\cal P}_K \cong {\cal P}_{K^\vee}$
and $\nu_{K^\vee}\cong \nu_K$. % and ${\cal P}_K^\vee \cong {\cal P}_K$.
\item $0\leq \nu_K \leq \dim(X)$, and $\nu_K=\dim(X)$ iff $K$ is translation invariant under $X$.
\item $\nu_K = 0$ iff $K$ is in $M(X)$, i.e. iff $\chi(K)\neq 0$ holds for the Euler characteristic. In this case ${\cal P}_K=\delta_0$.
\item $\nu_K > 0$ iff $K$ is in $E(X)$, i.e. iff $\chi(K)=0$.
\item $(P,\nu)=({\cal P}_K,\nu_K)$ is uniquely characterized by the property: 
$P[-\nu]$ is a summand of $K^\vee*K$ with ${\cal H}^{-\nu}(P)_0 \neq 0$. We remark that then this stalk is dual to $End_{\bf D}(K)$ and hence 
isomorphic to $\Lambda$ (see also [BN, cor.2]). 
%NO and ${\cal H}^{i}(K*K^\vee)=0$ for $i > -\nu$.  NO
\item ${\cal P}_K[\pm \nu_K]$ has multiplicity one in $K*K^\vee$ and ${\cal P}_K[n] \hookrightarrow K^\vee*K$ implies $\vert n\vert \leq \nu_K$.
\item $x\in supp({\cal H}^0({\cal P}_K[-\nu_K]))$ iff $T_x^*(K) \cong K$ (this describes the stabilizer of $K$). %\item  $\nu_K +\dim(supp({\cal P}_K)) \geq \dim(X)$.
\item $\nu_{K_\chi} = \nu_K $, $\mu(K)=\mu(K_\chi)$ and 
${\cal P}_{K_\chi} = ({\cal P}_{K})_\chi$ for twists $K_\chi$ of $K$.
\item $\nu_K = \mu({\cal P}_K)$.
\end{enumerate}
where for a complex $G$ in $\bf D$ we define  $$\mu(G)= \max\{\nu \ \vert \ {\cal H}^{-i}(G) = 0 \mbox{ for all } \ i < \nu \}\ .$$
\end{Lemma}

\medskip
{\it Proof}. For property 8 use that twisting with a character defines a tensor functor (see [KrW]).
The symmetry isomorphism $S: K^\vee * K \cong K * K^\vee$ together with property 5
gave ${\cal P}_{K^\vee}\cong {\cal P}_K$ and $\nu_{K^\vee}=\nu_K$. For property 2
notice that the perverse cohomology of 
the direct image $Ra_*(K^\vee\boxtimes K)$ vanishes 
in degree $> \dim(X)$, and for $\nu_K=\dim(X)$ one easily shows $a^*({\cal P}_K)[\dim(X)]\cong K^\vee\boxtimes  K$. 
%and $X\times \{ 0\}$ 
Hence  ${\cal P}_K[\dim(X)] \cong K \otimes {\cal H}^\bullet(K^\vee)_0$ by restriction to $\{0\} \times X$, and hence ${\cal P}_K\cong K$. Then $T_x^*(K)\cong K$ for all $x\in X$ follows by restriction to $\{x\} \times X$. 
The proof of property 3 and 4  follows from the next commutative diagram, whose right side stems from the hard Lefschetz theorem (see also [BN, 2.6])
$$ \xymatrix@+0,3cm{    \delta_0   &     {\cal P}_K[-\nu_K]  \ar@<1ex>@{^{(}->}[d]\ar[l]_{\exists !\ \epsilon} \cr
K^\vee*K \ar[u]^{eval_K}  &      \bigoplus_{\nu =0}^{\nu_K}  {\cal P}_K[2\nu-\nu_K]\ \ \oplus\ rest \ar[l]_-\sim \ar@<1ex>@{->>}[u]^{pr_{-\nu_K}} \cr
K*K^\vee \ar[u]^{S} \ar[r]^-\sim &      \bigoplus_{\nu = 0}^{\nu_K}  {\cal P}_K[2\nu-\nu_K] \ \oplus\ \ rest \ar[u]_{\sim} \ar@<1ex>@{.>>}[d]\cr
\delta_0 \ar@/^23mm/[uuu]^{\chi(K)}  \ar[u]^{coev_K}\ar[r]^{\exists !\ \sigma}  &     {\cal P}_K^\vee[+\nu_K]  \ar@<1ex>@{^{(}.>}[u]^{i_{+\nu_K}}}
$$
The two middle horizontal arrows  define $ \bigoplus_{\nu = 0}^{\nu_K}  {\cal P}_K[2\nu-\nu_K] $ as a retract of $K*K^\vee$, using $S:K*K^\vee \cong K^\vee*K$.  The middle vertical arrow on the right is an isomorphism respecting the direct sum decomposition
 $\bigoplus_{\nu = 0}^{\nu_K}  {\cal P}_K[2\nu-\nu_K]\ \oplus\ rest$. 
The existence of such a decomposition follows from the hard Lefschetz theorem, since the symmetry $S$ can be chosen so that it commutes with the Lefschetz maps $L$. Indeed, by defining the Lefschetz morphism $L: K^\vee *K \to K^\vee*K[2](1)$ as 
$L=Ra_*(\eta)$, where $\eta: K^\vee \boxtimes K \to K^\vee \boxtimes K[2](1))$ is induced by 
the cup-product of $K^\vee \boxtimes K$ with the morphism $\Lambda \to \Lambda[2](1)$
defined by an ample theta divisor of $X\times X$ whose Chern class is symmetric with respect to
the switch $\sigma_{12}(x_1,x_2)=(x_2,x_1)$, it suffices to know that $S=Ra_*(\psi)$ holds for some
isomorphism $\psi: K^\vee \boxtimes K \cong \sigma_{12}^*(K \boxtimes K^\vee)$. For this
see [BN, 2.1]. Since $coev_K: \delta_0 \to K^\vee*K$ factorizes over ${\cal P}_K^\vee[+\nu_K]$
and since $Hom_{\bf D}({\cal P}_K^\vee[+\nu_K], {\cal P}_K[-\nu_K])$ vanishes unless $\nu_K \leq - \nu_K$ and hence
$\nu_K=0$,
this proves assertion 3 and 4 taking into account the discussion of the case $\nu_K=0$
given in [KrW]. 
\qed

\medskip
We will show  ${\cal P}_K^\vee \cong {\cal P}_K$ later in lemma \ref{17}. Using this already, 
the lower right part of the last diagram is contained in $\bigoplus_{\nu = 0}^{\nu_K}  {\cal P}_K[2\nu-\nu_K]$ using the fact that both ${\cal P}_K[\pm \nu_K]$ appear with multiplicity one
as a direct summand in $K^\vee*K$. Notice, that both morphisms $\epsilon$ and $\sigma$ are nontrivial morphisms in the category $\bf D$.

\medskip
Besides the above large  \lq{monoidal component}\rq\ diagram there are similar commutative diagrams for semisimple perverse objects $P$ in ${\bf P}$.

\medskip
For $P=\bigoplus_i m_iP_i$ and irreducible $P_i\in {\bf P}$ such that $P_i\not\cong P_j$ for $i\neq j$ there are commutative diagrams
$$ \xymatrix{   \delta_0 \ar@{=}[r] & \delta_0   &    \ \bigoplus_i m_i^2\cdot {\cal P}_{P_i}[-\nu_{P_i}]  \ar[l]_-{\sum_i\ tr\circ \epsilon_i} \cr
P^\vee*P \ar[u]^{eval_P}  & \bigoplus_i \ m_i^2 \cdot  P_i^\vee*P_i \ar@{_{(}->}[l]\ar[u]_{\sum_i eval_{P_i}}  &     \ \bigoplus_i \bigoplus_{\nu = 0}^{\nu_{P_i}}  m_i^2 \cdot  {\cal P}_{P_i}[2\nu-\nu_{P_i}]  \ar@{_{(}->}[l]\ar@{->>}[u]_{\bigoplus_i pr_{-\nu_{P_i}}}    }
$$

\medskip
Also the following diagrams are commutative.  Notice, part of the next diagram is displayed already  in the last diagram. However, the next two diagrams are commutative also in the reverse direction, i.e. with the additional arrows inserted. This follows from $Hom_{\bf D}(P_i^\vee*P_j,\delta_0)
= Hom_{\bf D}(P_j,P_i)=0$ for irreducible $P_i \not\cong P_j$ in $\bf P$. 
The lower diagram is obtained from the upper one 
by Tannaka duality

$$ \xymatrix{   \delta_0 \ar@{=}[r] & \delta_0   &    \ \bigoplus_i m_i^2\cdot  {\cal P}_{P_i}[-\nu_{P_i}]  \ar[l]_-{\sum_i\ tr\circ \epsilon_i} \ar@<-1ex>@{_{(}->}[d]\cr
P^\vee*P \ar@<1ex>@{->>}[r]\ar[u]^{eval_P}  & \bigoplus_i \ m_i^2 \cdot  P_i^\vee*P_i \ar@{^{(}->}[l]\ar@<1ex>@{->>}[r]\ar[u]_{\sum_i eval_{P_i}}  &     \ \bigoplus_i \bigoplus_{\nu = 0}^{\nu_{P_i}}  m_i^2\cdot  {\cal P}_{P_i}[2\nu-\nu_{P_i}]  \ar@{^{(}->}[l]\ar@{->>}[u]_{\bigoplus_i pr_{-\nu_{P_i}}} \cr
P*P^\vee  \ar@{->>}[r] & \bigoplus_i \ m_i^2\cdot  P_i*P_i^\vee \ar@{->>}[r] \ar@<1ex>@{^{(}->}[l]&    \  \bigoplus_i\bigoplus_{\nu = 0}^{\nu_{P_i}} \ m_i^2\cdot  {\cal P}_{P_i}^\vee[2\nu-\nu_{P_i}] \ar@<1ex>@{^{(}->}[l]\ar@<-1ex>@{->>}[d] \cr
  \delta_0 \ar@{_{(}->}[u]^{coev_P} \ar@{=}[r] & \delta_0 \ar[u]_{\oplus_i coev_{P_i}}\ar[r]^-{\oplus_i\ tr^\vee \circ \sigma_i}  &  \ \bigoplus_i  \ m_i^2\cdot {\cal P}_{P_i}^\vee[+\nu_{P_i}]  \ar@{_{(}->}[u]_{\bigoplus_i i_{+\nu_{P_i}}}   }
$$

\medskip
{\bf Example}. For irreducible perverse sheaves $K$ on $A$ and $L$ on $B$ and $K\boxtimes L$ on $A\times B$, we have
${\cal P}_{K\boxtimes L} = {\cal P}_K
\boxtimes {\cal P}_L$ and $\nu_{K\boxtimes L} = \nu_K + \nu_L$
so that ${\cal H}^{-i}({\cal P}_K\boxtimes {\cal P}_L)$    is a skyscraper sheaf with stalk cohomology ${\cal H}^{-\nu_K-\nu_L}({\cal P}_K\boxtimes {\cal P}_L)_0
\cong \Lambda$ at the point zero for $i=\nu_K+\nu_L$, and vanishes for $i< \nu_K + \nu_L$.  
%In particular $K\boxtimes L \in F(A\times B)$ for
%$K\in F(A)$ and $B\in F(B)$.  

\medskip
{\bf Example}. An irreducible perverse sheaf $K$ is {\it negligible}\footnote{An equivalent definition is, that there exists an isogeny $g: A\times B \to X$ such that
$g^*(K)= \tilde K \boxtimes \Lambda_B[\dim(B)] $ for some abelian subvariety $B\neq 0$ and some $\tilde K\in Perv(A,\Lambda)$} if it has the form $$K\cong \delta_B^\psi * M \quad , \quad \delta_B^\psi := i_*(\Lambda_B)[\dim(B)]_\psi $$ for an irreducible $M\in M(X)$ (see also [KrW]), a nontrivial abelian subvariety $i:B\hookrightarrow X$ and a twist 
by a character $\psi: \pi_1(X,0)\to \Lambda^*$. Then $K*K^\vee \cong  (H^\bullet(X,\delta_B)\otimes_\Lambda\delta_B^\psi) * (M*M^\vee)$. This allows to compute $coev_B$ and $coev_M$ separately. Hence,  the monoidal component is $$ \fbox{$ {\cal P}_K \cong \delta_B^\psi \quad , \quad \nu_K = \dim(B) $} \ $$  by assertion 3) of lemma \ref{1}. 
 Indeed for an irreducible perverse sheaf $M\in M(X)$ we have ${\cal P}_M =\delta_0$.
The above formula for $\nu_K$  is a special case of  $$ \nu_F= \nu_K + \dim(A) \quad , \quad F=p^*(K)[\dim(A)]$$ for quotient morphisms $p:X\to B=X/A$, which  by an isogeny is easily reduced  to the case $X=A\times B$ where $p$ is the projection to the second factor and  $F=\delta_A \boxtimes K$. Then the assertion is obvious. Indeed, for $K\boxtimes L$ on $A\times B$ and $K\in Perv(A,\Lambda)$ and $L\in Perv(B,\Lambda)$, one has $\nu_{K\boxtimes L}  = \nu_K  + \nu_L$.

\medskip 
{\bf Tensor ideals}.
Semisimple complexes, whose irreducible perverse constituents (with shifts)
are translation invariant by nontrivial abelian subvarieties, resp. whose constituents
have Euler characteristic zero, define tensor ideals ${\bf N}$ and ${\bf N}_E$
in the tensor category $({\bf D},*)$ so that ${\bf N} \subset {\bf N}_E$. One can show that a complex $K$ is translation invariant under an abelian subvariety $A\subseteq X$ iff all perverse constituents of all 
perverse cohomology sheaves ${}^p H^i(K)$ are translation invariant under $A$. Furthermore by lemma \ref{1}, assertion 7 for $\nu_K>0$ an irreducible perverse sheaf is in ${\bf N}_E$ but not in ${\bf N}$ iff
${\cal H}^0({\cal P}_K[-\nu_K])$ is a skyscraper sheaf. Let $E(X)$ resp. $N(X)$
denote the perverse sheaves in ${\bf N}_E$ resp. ${\bf N}$, and $F(X)$ the isomorphism classes of irreducible perverse $K$ in $E(X)\setminus N(X)$.

\medskip
{\bf Reconstruction}.  We know
$Hom_{\bf D}(K*K^\vee, {\cal P}_K[-\nu_K])\cong \Lambda\neq 0$ for irreducible $K$ in ${\bf P}$. 
By rigidity $Hom_{\bf D}(K, K*{\cal P}_K[-\nu_K])\cong \Lambda \neq 0$, so there exists
a nontrivial morphism $K[\nu_K]\to  {\cal P}_K *K \cong K*{\cal P}_K$. Our aim is to show that there exists
a retract morphism  in ${\bf D}$ (of course unique up to a scalar)
$$  K \hookrightarrow {\cal P}_K[-\nu_K] *K \ .$$ 
Similarly, by rigidity then ${\cal P}_K^\vee*K \neq 0$.
By the decomposition theorem ${\cal P}_K^\vee*K = \bigoplus_{L,\nu}
L[-\nu]$ decomposes into a sum of shifted irreducible perverse sheaves $L$ (with $\nu\in\Z$). By the rigidity and strictness of the additive category ${\bf D}$
the morphism $id_K: K=\delta_0*K \to K*K^\vee*K \to K*\delta_0=K$ \lq{factorizes}\rq\  
in the form $id_K = \sum_{L,\nu} {\rm v}_{L,\nu} \circ u_{L,\nu}$.
The left horizontal morphism in nthe next diagram is the composite of $\varphi = \sigma * id_K$ and the monomorphism $\iota' * id_K$
$$ \xymatrix@+0,5cm{ K= \delta_0 * K\ar[r]^{coev_K * id_K}\ar[rd]_{\varphi}  \ar@{.>}[rdd]_{\exists \ u} &  K*K^\vee*K \ar[r]^{id_K * eval_K} & K *\delta_0 =K \cr
& {\cal P}_K^\vee[\nu_K] * K \ar@{^{(}->}[u]  \ar[ru]_\psi  &  \cr
& L[\nu_K-\nu] \ar@{^{(}->}[u] \ar[ruu]_{\rm v} & } $$
Then $\psi \circ \varphi= id_K$ for $\psi= (id_K *eval_K) \circ (\iota'*id_K)$.
Therefore $\varphi\neq 0$, and for  some constituent $i: L[-\nu]\hookrightarrow {\cal P}_K^\vee*K$  there exists a nontrivial morphism $u=u_{L,\nu}$ so that for ${\rm v}={\rm v}_{L,\nu}$ in ${\bf D}$ as in the diagram ${\rm v}\circ u\neq 0$. Warning: Notice $v= \psi\circ i$, but the lower left of the diagram may not be commutative. If $r$ is a retract of $i$, then $u=r\circ \varphi$.

\medskip
%$L[-\nu+\nu_K]\hookrightarrow {\cal P}_K[\nu_K]$, in the sense that there are
%nontrivial maps u,v as in the diagram. 
Nontrivial morphisms v from $L[\nu_K-\nu]$ to $K*\delta_0$ 
in ${\bf D}$ exist only for  $\nu_K-\nu\leq 0$, nontrivial morphisms $u$ in ${\bf D}$ from $\delta_0*K$ to $L[\nu_K-\nu]$ only for
$\nu_K-\nu\geq 0$. Hence $\nu=\nu_K$. That ${\rm v}\circ u: K=\delta_0*K\to L \to K*\delta_0$ is nontrivial 
forces $u$ and v to be isomorphisms of perverse sheaves $L\cong K$, since both $L$ and $K$ are irreducible. This proves $K \hookrightarrow {\cal P}_K^\vee [\nu_K] *K $ or 
$K[-\nu_K] \hookrightarrow {\cal P}_K^\vee *K $, and hence by the hard Lefschetz theorem $$ \fbox{$ K[\pm \nu_K]\hookrightarrow {\cal P}_K^\vee *K $} \ .$$ 
Applying this for $K^\vee$ instead of $K$,  by passing to the
 Tannaka duals  we then obtain from lemma \ref{1}, part 1 the desired assertion  $$ \fbox{$ K[\pm\nu_K] \hookrightarrow {\cal P}_K *K $} \ .$$ 
 Together with $$ \fbox{$ {\cal P}_K[\pm \nu_K] \hookrightarrow K*K^\vee $} \ $$ this implies 
 
\begin{Lemma} For irreducible perverse sheaves $K$ in ${\bf P}$ and an abelian subvariety $A\subseteq X$ and homomorphisms $f:X\to Y$ the following holds 
\begin{enumerate}
\item $Rf_*(K)=0$ iff $\ Rf_*({\cal P}_K)=0$. 
\item $K\in E(X)$ iff ${\cal P}_K\in E(X)$.  
\item $K$ is invariant under $A\subseteq X$ iff ${\cal P}_K$ is invariant under $A$.
\item $K\in N(X)$ iff ${\cal P}_K\in N(X)$. 
\item $K\in F(X)$ iff ${\cal P}_K \in F(X)$.
\end{enumerate}
\end{Lemma}

{\it Proof}. Obviously $3) \Longrightarrow 4)$
and $2),4) \Longrightarrow 5)$. For 1) use that $Rf_*$ is a tensor functor, for 2) use the hereditary
property of the class $N_{Euler}$ (see [KrW]),
and for 3) use $T^*_x(K*L)= T_x^*(K)*L$ for translations $T_x(y)=x+y$ for closed points $x\in A$
together with  $K[\pm\nu_K] \hookrightarrow {\cal P}_K *K$ and ${\cal P}_K[\pm \nu_K] \hookrightarrow K*K^\vee$. \qed

\medskip
{\bf Extremal perverse sheaves}. For closed points $x\in X$
the skyscraper sheaves $\delta_x$ are in ${\bf P}$ and $T_x^*(\delta_x)=\delta_0$, and $K\in {\bf P}$ iff $T^*_x(K)\in {\bf P}$. 
For $K,L$ in ${\bf P}$ the $\Lambda$-dual of the stalk cohomology $ {\cal H}^n(L^\vee*K)_x$ at  $x$ can be identified with $Hom_{\bf D}(L^\vee*K,\delta_x[-n])  \cong Hom_{\bf D}(T^*_x(K),L[-n])$, which is zero for $n>0$ by the perverse vanishing conditions for morphisms. Hence
${\cal H}^{>0}(K*{\cal P}_K) =0$. Since
$K[\pm \nu_K] \hookrightarrow K*{\cal P}_K$, therefore
${\cal H}^{>-\nu_K}(K) = {\cal H}^{>0}(K[-\nu_K])$ vanishes. 
For irreducible perverse sheaves $K$ this implies the inequality
$$  \fbox{$ \nu_K \ \leq\ \mu(K) $} \ .$$

%\medskip
%{\bf Remark}. In the case $char(k)=0$ for {minimal} $K\in F(X)$ the 
%stronger inequalities  $  \nu_K\leq \mu(X) \leq\mu(K) $ holds. See [W].

\medskip
Suppose $K\in F(X)$ is {\it extremal} in the sense that ${\cal H}^0(K[-\nu_K])\neq 0$,  or equivalently that $\nu_K =\mu(K)$, holds. We claim that
$$ \fbox{$ \nu_K=\mu(K) \ \ \Longleftrightarrow \ \ T^*_x(K) \cong {\cal P}_K \ \mbox{ for some } x\in X$} \ .$$
The implication $\Leftarrow$  follows from lemma 1, part 9. For the converse recall that $K[- \nu_K] \hookrightarrow {\cal P}_K*K$ and also $K[- \nu_K] \hookrightarrow {\cal P}_K^\vee *K$. 
Therefore  $\nu_K=\mu(K)$ implies  $0\neq  {\cal H}^0(K[-\nu_K]) $, and hence ${\cal H}^0({\cal P}_K^\vee*K)\neq 0$.
Notice, both $K$ and $L={\cal P}_K^\vee$ are irreducible perverse sheaves
and for irreducible perverse sheaves $K$ and $L$ one has ${\cal H}^0(L*K)\neq 0$
iff $T^*_x(K)\cong L^\vee$ holds for some $x\in X$ (see [BN,2.5], or the computations above). This implies $T^*_x(K)\cong {\cal P}_K$ for some $x\in X$, and proves our claim.

\medskip
If $K={\cal P}$ is a monoidal perverse sheaf, then $K$ is extremal and furthermore ${\cal H}^{-\nu_K}(K)_0\neq 0$ holds. Therefore the argument above shows that we even get an isomorphism $K\cong {\cal P}_K$, indeed we get this for $x=0$ from
the stronger assertion ${\cal H}(K[-\nu_K])_0 \subseteq {\cal H}^0({\cal P}_K*K)_0$. 
Furthermore, the same argument then
applied to the retract $K[- \nu_K] \hookrightarrow {\cal P}_K*K$, instead of  $K[- \nu_K] \hookrightarrow {\cal P}_K^\vee*K$, shows $K \cong {\cal P}_K^\vee$. Therefore $K^\vee \cong K$ follows for monoids $K$.
 
\medskip
Using this information, we get ${\cal P}_K^\vee \cong {\cal P}_K$ for arbitrary irreducible $K\in {\bf P}$. Hence if $K$ is extremal, then $K^\vee \cong T^*_{2x}(K)$. If $K$ is extremal and self dual in the sense $K^\vee \cong K$, then 
$T^*_{2x}(K)\cong K$. If $K$ is a monoidal component, then $K$ is extremal. Altogether this implies

\begin{Lemma} \label{17} For an irreducible perverse sheaf 
$K$ one has $\nu_K\leq \mu(K)$. $K$ is extremal in the sense $\nu_K=\mu(K)$ iff $K$ is isomorphic to a translate of its monoidal component. If $K$ is the monoidal component of an irreducible perverse sheaf, then  $$ \fbox{$ K^\vee \cong K\cong {\cal P}_K$} \ .$$
In particular, we obtain $\nu_K = \mu({\cal P}_K) = \nu_{{\cal P}_K}$.
\end{Lemma}

For monoids $K={\cal P}={\cal P}_K$ we have the following
commutative diagram, using that $K[\pm d] \hookrightarrow K*K$
occurs with multiplicity one in $K*K$ and also using $K^\vee \cong K$
$$ \xymatrix@+0,8cm{
K \ar@{.>}[ddr]_{\sim}\ar[r]^-{coev_K * id_K} \ar[dr]^{\sigma*id_K} &  (K*K^\vee)*K \ar[r]^-{ass}_\sim & K*(K^\vee*K) \ar@{->>}[d]_{id_K*j^\vee} \ar[r]^-{id_K*eval_K} & K \cr
 &  K[+\nu_K]*K \ar[r]^{\exists ! \ a}\ar@{^{(}->}[u]_{j*id_K}   & K*K[-\nu_K] \ar[ru]^{id_K*\epsilon}\ar@{->>}[d]_{p'}  \cr
 & K \ar@{^{(}->}[u]_{j'} \ar[r]^\sim & K \ar@{.>}[uur]_{\sim}& \cr } $$
for the diagrams \label{smalldiagrams}
$$ \xymatrix{   \delta_0 \ar[dr]_{\sigma} \ar[r]^-{coev_K} & K*K^\vee \cr
& K[\nu_K] \ar@{^{(}->}[u]_j } \quad \quad \quad \xymatrix{    K^\vee*K \ar@{->>}[d]_{j^\vee} \ar[r]^-{eval_K} & \delta_0   \cr
 K[-\nu_K] \ar[ur]_{\epsilon} & }$$
The two small diagrams, together with rigidity, imply the existence of $a$ such that $id_K = (id_K*\epsilon) \circ a \circ (\sigma*id_K)$. Repeating the argument, used in the section on reconstruction, for $\varphi=\sigma*id_K$ and 
$\psi=(id_K*\epsilon)\circ a$, we see that $id_K= \psi\circ \mu\circ (\sigma*id_K)$ factorizes over the unique (!) retract $\mu: K[+\nu_K]*K \twoheadrightarrow K$  to the unique constituent $j':K\hookrightarrow K[+\nu_K]*K$ isomorphic to $K$. Similarly, there is a unique retract $p': K*K[-\nu_K] \twoheadrightarrow K$. Repeating the argument, used in the section on reconstruction, now for $\varphi=a\circ (\sigma*id_K)$ and $\psi=id_K*\epsilon$, we find a commutative diagram \label{funnydiag}
$$ \xymatrix@+0,8cm{
K \ar[ddr]_{\sim}\ar[r]^-{coev_K * id_K} \ar[dr]^{\sigma*id_K} &  (K*K^\vee)*K \ar[r]^-{ass}_\sim & K*(K^\vee*K) \ar@{->>}[d]_{id_K*j^\vee} \ar[r]^-{id_K*eval_K} & K \cr
 &  K[+\nu_K]*K \ar@{->>}[d]^\mu \ar@{^{(}->}[u]_{j*id_K}   & K*K[-\nu_K] \ar[ru]^{id_K*\epsilon}  \cr
 & K  \ar[r]^\sim & K \ar@{^{(}->}[u]^r\ar[uur]_{\sim}& \cr } $$
$\mu\circ (\sigma*id_K): K\to  K$ completes the left lower part of the diagram.

\medskip
{\bf Functors}. For $\Lambda$-linear  tensor functors $F$ between rigid  symmetric monoidal  (not necessarily abelian) $\Lambda$-linear  tensor categories,  $F(coev_K)=coev_{F(K)}$ and also $F(eval_{K})=eval_{F(K)}$ holds.
We will use this for the direct image functor $F=Rf_*$
which for a homomorphism  between abelian varieties 
$$f:X\to Y$$
induces a triangulated tensor functor between ${\bf D}(X) \subseteq D_c^b(X,\Lambda)$
and ${\bf D}(Y) \subseteq D_c^b(Y,\Lambda)$.

\smallskip
{\bf Assumption}. Suppose 
$Rf_*(K)$ is perverse\footnote{For what follows one also could replace $\bf D$ by  some localization ${\bf D_H}$ with respect to a hereditary class ${\bf H}$ (see [KrW]), and then it suffices to assume $Rf_*(K)\in {\bf P}_H$. For complex abelian varieties on the other hand the assumption can always be achieved by a generic character twist using the relative vanishing theorem of [KrW].}. By the decomposition theorem $Rf_*(K)$ decomposes into nonisomorphic irreducible perverse sheaves $P_i$ with multiplicities $m_i$
$$ \fbox{$ L:= Rf_*(K) = \bigoplus_i \ m_i\cdot P_i $} \ .$$
Since $Rf_*$ is a tensor functor
$$ Rf_*(K*K^\vee) \ = \ L*L^\vee \ =\  \bigoplus_i \  m_i^2 \cdot P_i*P^\vee_i \ \  \oplus\
\ \bigoplus_{i\neq j}\ m_im_j\cdot P_i*P^\vee_j \ .$$
Using lemma \ref{1}, property 5 of monoidal components and the adjunction formulas from page \pageref{adformeln} it is easy to see that any irreducible constituent $Q$ of $L*L^{\vee}$ with
${\cal H}^0(Q)_0\neq 0$ is contained in the first sum 
$\bigoplus_i  m_i^2 \cdot P_i*P^\vee_i$, hence is of the form 
$$Q\ \cong\ {\cal P}_{P_i}[-\nu_{P_i}] \ .$$
Now applying $F=Rf_*$ to the monoidal diagram of $K$ 
gives the right side of the following commutative diagram
$$ \xymatrix@+0,3cm{ \bigoplus_i m_i^2\cdot {\cal P}_{P_i}[-\nu_{P_i}] 
\ar[r]^-{\oplus_i\ tr\circ \epsilon_i}  
& \delta_0   &     Rf_*{\cal P}_K[-\nu_K]  \ar[l]_{Rf_*(\epsilon)} \cr
\bigoplus_i \bigoplus_{\nu =0}^{\nu_{P_i}}  m_i^2\cdot {\cal P}_{P_i}[2\nu -\nu_{P_i}]  
\ar@{^{(}->}[r]\ar@{->>}[u]_{pr}  
  & L^\vee*L \ar[u]^{eval_P}  
 &      \bigoplus_{\nu = 0}^{+\nu_K}  Rf_*{\cal P}_K[2\nu-\nu_K]\ \ \oplus\ rest \ar[l]_-\sim\ar@{->>}[u]_{Rf_*(pr_{-\nu_K})} \cr
\bigoplus_i\bigoplus_{\nu = 0}^{\nu_{P_i}} \ m_i^2 \cdot {\cal P}_{P_i}[2\nu-\nu_{P_i}]   \ar@<1ex>@{^{(}->}[r]
 &  L*L^\vee \ar[u]^{S}_\sim \ar[r]^-\sim  \ar@{->>}[l] 
 &      \bigoplus_{\nu = 0}^{\nu_K} Rf_* {\cal P}_K[2\nu-\nu_K] \ \oplus\ \ rest \ar@/_12mm/[l]^\sim\ar[u]_{\sim} \cr 
 \bigoplus_i  \ m_i^2 \cdot {\cal P}_{P_i}[+\nu_{P_i}]  \ar@{^{(}->}[u]_{i}  
& \delta_0 \ar[u]^-{coev_P}\ar[r]^{Rf_*(\sigma)} \ar[l]_-{\oplus_i\ tr^\vee \circ \sigma_i}
&    Rf_* {\cal P}_K[+\nu_K]  \ar@{^{(}->}[u]_{Rf_*(i_{+\nu_K})} \cr    }$$

\medskip
The lower part of this diagram defines the next commutative diagram
$$ \xymatrix@+0,5cm{ 
\bigoplus_i\bigoplus_{\nu =0}^{\nu_{P_i}} \ m_i^2\cdot {\cal P}_{P_i}[2\nu -\nu_{P_i}]  \ar@{->>}@/_8mm/[d]_{p}  \ar@<-1ex>@{_{(}->}[r]_-{\rm v'}
 &  L*L^\vee\ \ar@{->>}[l]_-{\rm v}  \ar@<-1ex>@{->>}[r]_-{u'}  &      \ \bigoplus_{\nu = 0}^{\nu_K} Rf_* {\cal P}_K[2\nu-\nu_K] \ar@{_{(}->}[l]_-{u}  \ar@{->>}@/_8mm/[d]_{Rf_*(pr)}
  \cr 
 \bigoplus_i  \ m_i^2 \cdot {\cal P}_{P_i}[+\nu_{P_i}]  \ar@{^{(}->}[u]_{i}  
& \delta_0 \ar[u]^{coev_P}\ar[r]^{Rf_*(\sigma)} \ar[l]_-{\oplus_i\ tr^\vee \circ \sigma_i}&    Rf_* {\cal P}_K[+\nu_K]  \ar@{^{(}->}[u]_{Rf_*(i_{+\nu_K})} \cr    }$$
where the retract morphism $u$ is obtained from the right middle diagram, using the isomorphism $S$ and taking into account 
that the coevaluation map of $L$ ignores the part of the last diagram
entitled \lq{rest}\rq. Altogether this defines a morphism
$$ Rf_*{\cal P}_K[+\nu_K] \ \longrightarrow \ \bigoplus_i\bigoplus_{\nu =0}^{\nu_{P_i}} \ m_i^2\cdot {\cal P}_{P_i}[2\nu -\nu_{P_i}] $$
whose \lq{image}\rq\ is contained in $\bigoplus_i  \ m_i^2 {\cal P}_{P_i}[+\nu_{P_i}] $ and, without loosing information, can be considered as a morphism
$ Rf_*{\cal P}_K[+\nu_K] \ \longrightarrow \ \bigoplus_i \ m_i^2\cdot {\cal P}_{P_i}[+\nu_{P_i}]$.
For $L\neq 0$, from the definition it is clear   that for each $i$ the  composed morphism
$$  Rf_*({\cal P}_K)[+\nu_K] \ \longrightarrow \ \bigoplus_i \ m_i^2 \cdot{\cal P}_{P_i}[+\nu_{P_i}]\ \longrightarrow\ m_i^2 \cdot {\cal P}_{P_i}[+\nu_{P_i}] $$
is nontrivial. Indeed, if $pr_i\circ p \circ {\rm v} \circ u \circ i_{+\nu_K}$ would be zero, then 
also the composition with $Rf_*(\sigma)$, which is $tr^\vee \circ \sigma_i\neq 0$, would be zero.   
The same argument also implies $Rf_*(\sigma)\neq 0$. Hence we can repeat this argument in the other direction to show that the composed morphism $Rf_*(pr)\circ u' \circ {\rm v}' \circ i$
$$ m_i^2 \cdot {\cal P}_{P_i}[+\nu_{P_i}] \longrightarrow Rf_*({\cal P}_K[+\nu_K]) $$
is again nontrivial, and also their composition. 
This proves

\begin{Proposition} \label{myst} Suppose $K$ is  an irreducible perverse sheaf so that the semisimple complex $L=Rf_*(K)=\bigoplus_{i\in I} m_i\cdot P_i$ is perverse and not zero (i.e $m_i>0$).
Then for every irreducible perverse constituent $P_i$ of $L$
there exist nontrivial morphisms in the derived category $$Rf_*({\cal P}_K)[\nu_K] \ \ \longrightarrow\ \  m_i^2\cdot {\cal P}_{P_i}[\nu_{P_i}] \  $$ 
$$ m_i^2\cdot  {\cal P}_{P_i}[\nu_{P_i}]  \ \ \longrightarrow\ \    Rf_*({\cal P}_K)[\nu_K] \ $$
whose composition (in both directions) is not zero.
\end{Proposition}

\medskip
Then $Hom_{{\bf D}}(M,N[r]) =0$ for perverse sheaves $M,N$ and $r<0$  implies

\begin{Corollary} \label{absch} Suppose $K$ is an irreducible 
perverse sheaf for which the semisimple complexes $L=Rf_*(K)=\bigoplus_i m_i\cdot P_i$ and $Rf_*({\cal P}_K)$ are perverse with $L\neq 0$. Then $ \nu_K = \nu_{P_i}$ holds for all irreducible perverse constituents $P_i$ of $L$.
\end{Corollary}

Since $\nu_{P_i} \leq \dim(Y)$, we also obtain from proposition \ref{myst}
$$Rf_*(K)\neq 0 \mbox{ is in } {\bf P} \ \ \Longrightarrow \ \  \nu_K \leq \dim(Y) \ .$$

\medskip
{\bf Definition}. An irreducible perverse sheaf $F$ on $X$ will be called  {\it maximal},
if for {\it every} projection $f:X \to B$ to a simple quotient abelian variety $B$ of $X$ the direct images
$Rf_*(K_\chi)$ and $Rf_*({\cal P}_{K_\chi})$ are perverse and not zero for generic character twists $\chi$.
If $X$ is simple, any irreducible perverse sheaf $F$ is maximal. 

\medskip
{\it Example}. Perverse sheaves in $M(X)$ are maximal. 

\medskip
Define $\mu(X)$ to be the minimum of the dimensions of the (nontrivial) simple 
abelian quotient varieties $B\neq 0$ of $X$.

\begin{Lemma} \label{numu}  Suppose $K$ is a {maximal} irreducible perverse sheaf. If $Rf_*(K_\chi)$ and $Rf_*{\cal P}_{K_\chi}$ are perverse for $f:X \to B$ and $\dim(B)=\mu(X)$, then $$ \fbox{$ \nu_K \ \leq \ \mu(X) $} \ .$$
In case that $char(k)=0$, this holds for any maximal perverse sheaf $K$.
\end{Lemma}

{\it Proof}.
$\nu_{K_\chi}$ only depends on $K$, but not on $\chi$ (property 8). It is shown in the relative vanishing theorem of [KrW],  that  for $k$ of characteristic zero  one can  always assume that $L=Rf_*(K_\chi)=\bigoplus_i m_i\cdot P_i$ and $Rf_*({\cal P}_{K_\chi}) = Rf_*\bigl({\cal P}_{K\chi}\bigr)$ are perverse by applying a twist with a suitable generic character $\chi: \pi_1(X,0)\to \Lambda^*$. If $K$ is minimal, we
can therefore always dispose over the arguments from above. \qed

\medskip
We remark that twists with characters $\chi': \pi_1(B,0)\to \Lambda^*$ have the following effect:  $L=\bigoplus_i P_i$ changes into $P_{\chi'} = \bigoplus_i (P_i)_{\chi'}$,  ${\cal P}_{K}$ and ${\cal P}_{P_i}$ change as well into their $\chi'$-twist. This implies, that the morphisms constructed above are {\it independent} from twists of $K$ with characters $\chi'$ of $\pi_1(B,0)$.

\medskip
{\bf Functors revisited}. Suppose given a homomorphism $f:X\to Y$ of abelian varieties and  semisimple perverse sheaves $K$
and $P$ (or more genertally complexes)
on $X$ and some integer $\nu$ (by abuse of notation we then again write $\nu=\nu_K$) together with a commutative diagram
$$ \xymatrix{ \ \ K^\vee*K \ \ \ar@<1ex>@{->>}[dd]^p \ar[rd]^{eval_K} &  \cr
& \ \ \delta_0 \cr
\ \ P[-\nu_K]\ \ \ar@<1ex>@{^{(}->}[uu]^\iota \ar[ru]_\epsilon & } $$
such that $p\circ \iota=id$. Then $L=Rf_*(K)=\bigoplus_{i\in I} P_i[\lambda_i] $ and $Q=Rf_*(P)=
\bigoplus_{j\in J} Q_j[\lambda_j]$ decomposes with simple perverse sheaves $P_i$ and $Q_j$.
By abuse of notation, the index index sets $I$ and $J$ are not correlated to each other, so the same holds for the $\lambda_i$ and $\lambda_j$.  
With these notations 
we get 

\begin{Theorem} \label{ungl}
For any (shifted perverse) constituent $P_i[\lambda_i] \hookrightarrow L$ there exists
a (shifted perverse) constituent $Q_j[\lambda_j] \hookrightarrow Q$ such that 
$$ \nu_{Q_j} \leq \mu(Q_j) \leq \nu_K - \lambda_j \leq \nu_{P_i} \ $$
holds, and a constituent $P_{i'}[\lambda_{i'}] \hookrightarrow L$
such that $\nu_{P_{i'}}\leq \nu_K - \lambda_j$ holds.
%Furthermore $Q_j \cong {\cal P}_{P_i}$ holds.
%Probably  $\nu_P - \lambda_j = \nu_{P_i}$
\end{Theorem}

{\it Proof}. Since $Rf_*$ is a tensor functor, we get the commutative diagram
$$ \xymatrix{ \ \ L^\vee*L \ \ \ar@<1ex>@{->>}[dd] \ar[rd]^{eval_L} &  \cr
& \ \ \delta_0 \cr
\ \ Q[-\nu_K]\ \ \ar@<1ex>@{^{(}->}[uu] \ar[ru]_{Rf_*(\epsilon)} & } $$
for the evaluation morphism $eval_L$. For any direct factor $C=P_i[\lambda_i]$
in $L$ the evaluation morphism of $eval_C: C^\vee * C = P_i^\vee[-\lambda_i]*P_i[\lambda_i]=P_i^\vee * P_i \to \delta_0$ is induced by the evaluation morphism $eval_{P_i}$, which is computed 
via the upper horizontal morphisms of the next commutative diagram. 
The evaluation $eval_C$  is also
obtained as the restriction of the evaluation morphism $eval_L: L^\vee*L \to \delta_0$
to $C^\vee * C \hookrightarrow L^\vee*L$. The evaluation morphism $eval_L$ is given by the lower horizontal
morphisms of the next diagram. Altogether, this implies the existence of a morphism $\varphi$ 
$$ \xymatrix{ {\cal P}_{P_i}[-\nu_{P_i}]\ar[dr]_{\epsilon_i}  \ar[rr]^\varphi & & Q[-\nu_K]\ar[dl]^{Rf_*(\epsilon)} \cr
& \delta_0 & } $$ 
making the following diagram commutative 
$$ \xymatrix@+0,5cm{  P_i^\vee*P_i = C^\vee*C   \ar@<1ex>@{->>}[r]\ar@{^{(}->}[d] &   {\cal P}_{P_i}[-\nu_{P_i}] \ar[d]^{\varphi} \ar[r]^{\epsilon_i}\ar@<1ex>@{^{(}->}[l] &   \delta_0  \cr   L^\vee *L \ar@<1ex>@{->>}[r]  & Q[-\nu_K] \ar@<1ex>@{^{(}->}[l]\ar@{->>}[r] & \delta_0 \ar@{=}[u] } $$
Now we can decompose $Q[-\nu_K]=\bigoplus_{j\in J} Q_j[\lambda_j - \nu_K]$ and accordingly decompose also the morphism $\varphi$, so that for at least
one $j\in J$ we get a commutative diagram
$$ \xymatrix{ {\cal P}_{P_i}[-\nu_{P_i}]  \ar[dr]_{b_j}\ar[rr]^{a_j} & &  Q_j[\lambda_j -\nu_K] \ar[dl]^{c_j}\cr 
& \delta_0 & } \ .$$
with a {\it nontrivial} morphism $b_j$, since the morphism $\epsilon_i= \sum_j b_j$ is not zero. 
Then of course also $a_j\neq 0$ and $c_j\neq 0$.
Now $a_j\neq 0$ implies $-\nu_{P_i} \leq \lambda_j - \nu_K$, and $c_j\neq 0$ implies
${\cal H}^0(Q_j[\lambda_j - \nu_K])_0 = {\cal H}^{\lambda_j - \nu_P}(Q_j)_0 \neq 0$. Hence
$ - \nu_{P_i} \leq \lambda_j - \nu_K \leq - \mu(Q_j) \leq - \nu_{Q_j} $. 

\medskip
Reversing the argument,  we can conversely construct a {\it nontrivial} morphism
$$ \xymatrix{ Q_j[\lambda_j -\nu_K]   \ar[dr]_{}\ar[rr]^{} & &  {\cal P}_{P_i'}[-\nu_{P_i'}]\ar[dl]^{}\cr 
& \delta_0 & } \ $$
for some constituent $P_{i'}[\lambda_{i'}] \hookrightarrow L$. 
%Since $\varphi$ is a composition of morphisms of degree zero, therefore $\varphi$ and then also $a_j$ are also morphisms of degree zero. 
%Hence
%$-\nu_{P_i} = \lambda_j - \nu_P$ and ${\cal P}_{P_i} \cong Q_j$, since $a_j$ must be an isomorphism.
\qed
 
\medskip 
{\bf Definition}. Define $\nu_Q = \min_j (\nu_{Q_j})$ for the decomposition $Q= \bigoplus_j Q_j[\lambda_j]$,
and similarly define $\nu_L = \min_i (\nu_{P_i})$ for the decomposition $L= \bigoplus_i P_i[\lambda_i]$.

\medskip{\bf Definition}.  If $\nu_{Q_j}=\nu_Q$ (respectively $\nu_{P_i}=\nu_L$) holds,
a constituent $Q_j[\lambda_j]$ of $Q$ (respectively $P_i[\lambda_i]$ of $P$) will be called {\it minimal}.

\medskip
If we apply the last theorem
for a minimal constituent $P_i[\lambda_i]\hookrightarrow  L$, then the properties
$\nu_K - \lambda_j \leq \nu_{P_i} = \nu_L $
and $\nu_L \leq \nu_{P_{i'}} \leq \nu_K - \lambda_j$ imply
%$\nu_L \leq \nu_{Q_j}$ implies 
%$$  \nu_L = \nu_{Q_j} = \mu_{Q_j} = \nu_P - \lambda_j = \nu_{P_i} = \nu_L \ .$$
$$  \nu_{P_{i'}} = \nu_K - \lambda_j = \nu_{P_i} = \nu_L \ .$$
%In particular $Q_j$ is also minimal and $Q_j[\lambda_j - \nu_P] = Q_j[-\nu_L]$. Furthermore, the nontrivial morphism
% $$ a_j: {\cal P}_{P_i}[-\nu_{P_i}] = {\cal P}_{P_i}[-\nu_L] \longrightarrow Q_j[\lambda_j - \nu_P] = Q_j[-\nu_L] $$
%is then an isomorphism of shifted perverse sheaves. Hence
%$ Q_j\ \cong \ {\cal P}_{P_i}$.
In particular, $P_{i'}$ is also minimal and $Q_j[\lambda_j - \nu_K] = Q_j[-\nu_L]$. Furthermore, the nontrivial morphisms
 $$ a_j: {\cal P}_{P_i}[-\nu_{P_i}] = {\cal P}_{P_i}[-\nu_L] \longrightarrow Q_j[\lambda_j - \nu_K] = Q_j[-\nu_L] $$
and similarly
$$ Q_j[\lambda_j - \nu_K] = Q_j[-\nu_L] \to {\cal P}_{P_{i'}}[-\nu_{P_{i'}}] = {\cal P}_{P_{i'}}[-\nu_L] $$ 
imply $Q_j \cong {\cal P}_{P_{i'}}$. By lemma \ref{17}, the degree of 
an irreducible perverse sheaf is the degree of its monoidal perverse sheaf, and
we conclude for the degrees $$ \fbox{$ \nu_{Q_j} = \nu_{P_{i'}} = \nu_L $}\ .$$ Hence 
there is also an isomorphism of shifted perverse sheaves
$ Q_j\ \cong \ {\cal P}_{P_i}$. 

%Now assume $K=P$ and  then $L=Q= \bigoplus_i P_i[\lambda_i] = \bigoplus_j Q_j[\lambda_j]$.

\begin{Corollary} \label{mini}
%Suppose $K=P$. 
Let $P$ be the monoid attached to $K$ and $f:X \to Y$ be a homomorphism.
 For every minimal $P_i[\lambda_i] \hookrightarrow L=Rf_*(K)$ (i.e. $\nu_{P_i}= \nu_L$) 
there exists a shifted monoidal constituent $Q_j[\nu_K - \nu_L]\hookrightarrow Q=Rf_*(P)$ with
$$ Q_j \cong {\cal P}_{P_i} \ .$$ 
In particular, $\nu_{Q_j} = \nu_{P_i} = \nu_L$ and $\nu_Q \leq \nu_L$.  
\end{Corollary}

\medskip
{\bf An Application}. Let ${\cal P}$ be a monoidal perverse sheaf on $X$. Then $K=P={\cal P}\boxtimes {\cal P}$ is a monoid on $X\times X$ of degree $\nu_K= 2\nu_{\cal P}$.  For the morphism
$a: X\times X \to X$ we get  $L=Ra_*(K)={\cal P}*{\cal P}$. Since 
${\cal P}[-\nu_{\cal P}] \hookrightarrow {\cal P}*{\cal P}$ by lemma \ref{17}, this implies {\bf (*)}
$$  \nu_L = \min_i \nu_{P_i}  \leq \nu_{\cal P} \ .$$
By theorem \ref{ungl} and corollary \ref{mini}, the minimal constituents $P_i[\lambda_i]$ of $L$ 
give rise to monoidal constituents $Q_j[\lambda_j]\hookrightarrow L$ with the property
$ \lambda_j = \nu_K - \nu_L = 2\nu_{\cal P} - \nu_L $. Then, by the inequality {\bf (*)}, in particular
$$ 0 \leq \nu_{\cal P} \leq \lambda_j   \ .$$

\medskip
{\it Hard Lefschetz}. $Q_j[\lambda_j] \hookrightarrow L$ implies
$Q_j[\lambda_j-2i]\hookrightarrow L$ for all $i=0,\cdots ,\lambda_j$.
For $i=\nu_{\cal P} \leq \lambda_j $, therefore $$  Q_j[\lambda_j - 2\nu_{\cal P}] = Q_j[-\nu_L] \hookrightarrow L \ .$$
Notice ${\cal H}^0(Q_j[-\nu_L])_0 \neq 0$, since $Q_j$ is a monoidal perverse sheaf
and $\nu_{Q_j} = \nu_{P_i} = \nu_L$. 

\medskip
By lemma 1, part 5 there is a unique (shifted perverse) constituent in $L={\cal P}*{\cal P}$
with the property ${\cal H}^0(Q_j[-\nu_L])_0 \neq 0$, namely ${\cal P}[-\nu_{\cal P}]$. Hence $Q_j[-\nu_L]\cong {\cal P}[-\nu_{\cal P}]$ or $Q_j[\lambda_j]\cong {\cal P}[+\nu_{\cal P}]$. So $\lambda_j = \nu_{\cal P}$, in particular
$\nu_L = 2\nu_{\cal P} - \lambda_j = \nu_{\cal P}$.

\medskip
This proves

\begin{Lemma}\label{lfuenf}
For a monoid ${\cal P}$ on $X$ we have $\nu_{{\cal P}*{\cal P}} = \nu_{\cal P}$.
All (shifted perverse) constituents $Q_j[\lambda_j] \hookrightarrow {\cal P}*{\cal P}$ attached to a minimal (shifted perverse) constituent
$P_i[\lambda_i] \hookrightarrow {\cal P}*{\cal P}$ are isomorphic to ${\cal P}[+\nu_{\cal P}] \cong {\cal P}_{P_i}[\lambda_j]$ and minimal.
\end{Lemma}
 
\medskip
For $X$ consider the irreducible monoidal  perverse
sheaves  ${\cal P}$ on  $X$  with the property $\nu_{\cal P} < \dim(X)$. Let $\nu_+(X)$ be the maximum of all such $\nu_{\cal P}$. If $\nu_{\cal P}=\nu_+(X)$ holds, we call ${\cal P}$
a {\it maximal} monoid on $X$. 

\begin{Corollary}\label{maxim}
For a maximal irreducible monoid ${\cal P}$ on a simple abelian variety $X$ with ${\cal P}*{\cal P} \cong \bigoplus_i P_i[\lambda_i]$ either ${\cal P}_{P_i} \cong {\cal P}$ holds, or
${\cal P}_{P_i} \cong \delta_X^{\varphi_i}$ for some character $\varphi_i$. 
\end{Corollary} 
 
{\it Proof}. For $L={\cal P}*{\cal P}$ we have shown $\nu_L=\nu_{{\cal P}*{\cal P}} = \nu_{\cal P}$. Hence, for maximal
${\cal P}$ there are no (shifted perverse) constituents in $L={\cal P}*{\cal P}$ with degree $\nu_{P_i} > 
\nu_L$ except for ${\cal P}_{P_i} \cong \delta_X^{\varphi_i}$ by lemma 1, part 2. Hence every $P_i$ is either translation-invariant under $X$, or $\nu_{P_i}=\nu_{\cal P}=\nu_L$ is minimal in $L$. So we apply  Lemma \ref{lfuenf}.  \qed

\begin{Corollary}\label{maincor}
For monoids ${\cal P}_1,{\cal P}_2$ with degrees $\nu_1\leq \nu_2$ on an abelian variety $X$ with ${\cal P}_1 \not\cong {\cal P}_2$ the convolution $L={\cal P}_1 * {\cal P}_2$ has minimal degree $\nu_L > (\nu_1 +\nu_2)/2$.
\end{Corollary} 

\medskip{\it Proof}. We apply corollary  \ref{mini}
for the group law $a:X\times X\to X$ and $K=P={\cal P}_1\boxtimes {\cal P}_2$ with ${\nu_K}=\nu_1 + \nu_2$
and $L=Q=a_*(K)={\cal P}_1*{\cal P}_2 = \bigoplus_{i\in I} P_i[\lambda_i]$.   
Assume our assertion does not hold, i.e. suppose $\nu_L \leq (\nu_1 +\nu_2)/2$. This implies $\nu_L \leq \nu_2$ {\bf (*)}.
By corollary \ref{mini}, for any constituent $P_i[\lambda_i], i\in I$ with $\nu_{P_i}=\nu_L=\min_{i\in I}{\nu_{P_i}}$ there is  a monoidal constituent $Q_j[\lambda_j]$ in $L$
so that $\lambda_j = \nu_K - \nu_L = \nu_1 + \nu_2 - \nu_L$. The inequality {\bf (*)} implies $\lambda_j \geq 0$. Hence,  $Q_j[\lambda_j] \hookrightarrow L$ and $\lambda_j \geq 0$, by the Hard Lefschetz Theorem, also implies $Q_j[\lambda_j - 2i] \hookrightarrow L$ for all $i=0,..,\lambda_j$. For $i:= \lambda_j$, this gives the following constituent of $L$:
$$ Q_j[-\lambda_j] = Q_j[\nu_L - \nu_1 - \nu_2] = Q_j[-\nu_L][2\nu_L - \nu_1 - \nu_2] \hookrightarrow L \ .$$
By corollary \ref{mini}  we know that  $Q_j $ is a monoid with $ \nu_{Q_j}=\nu_L $. So from the above we conclude
$$  Q_j[-\nu_{Q_j}][2\nu_L - \nu_1 - \nu_2] \hookrightarrow L \ .$$
${\cal H}^0(Q_j[-\nu_{Q_j}])$ is a skyscraper sheaf with nontrivial stalk at $0$  and ${\cal H}^0({\cal P}_1*{\cal P}_2)_0\!\neq\! 0$ if and only if  ${\cal P}_1 \cong {\cal P}_2^\vee$ by [BN];
furthermore ${\cal H}^a({\cal P}_1*{\cal P}_2)_0=0$ for $a>0$. Since ${\cal P}_2^\vee \cong {\cal P}_2$ (Lemma \ref{17}), by our assumptions ${\cal P}_1 \not\cong {\cal P}_2^\vee$. Hence $
2\nu_L  - \nu_1 - \nu_2 $ must be $>0$. A contradiction.
\qed

\begin{Corollary}\label{EXTR}
For maximal monoids ${\cal P}_1,{\cal P}_2$ on an abelian variety $X$ with ${\cal P}_1 \not\cong {\cal P}_2$ the convolution ${\cal P}_1 * {\cal P}_2$ is translation invariant.
\end{Corollary}

\medskip
Similarly one obtains

\begin{Corollary}
For irreducible perverse sheaves $K_1,K_2$ on an abelian variety $X$ with degrees $\nu_{K_1} = \nu_{K_2} = \nu_+(X)$,
all simple constituents $P_i[\lambda_i]$ of $K_1*K_2$ are either translation-invariant under $X$, or ${\cal P}_{P_i} \cong {\cal P}_{K_1} \cong {\cal P}_{K_2}$. 
\end{Corollary} 

Also

\begin{Corollary}
For maximal irreducible monoids ${\cal P}_1,{\cal P}_2$ on a simple abelian variety $X$ with ${\cal P}_1 \not\cong {\cal P}_2$ assume ${\cal P}_1 * {\cal P}_2\neq 0$. Then $\delta^{\psi}[2\nu_+(X) - \dim(X)] \hookrightarrow {\cal P}_1 * {\cal P}_2$ for some character $\psi$. 
\end{Corollary}

\medskip
{\bf Isogenies}. We now discuss the behaviour of monoids 
with respect to pullback and push forward under isogenies $f:X\to Y$.

\medskip

\begin{Corollary} \label{W} Suppose $K$ is an irreducible monoidal perverse sheaf on $X$ 
with finite stabilizer $H=\{ x\in X\ \vert \ T_x^*(K)\cong K\}$. Then for the isogeny $\pi: X\to X/H$ 
 the direct image $L=\pi_*(K)$ is $L\cong \bigoplus_{\chi \in H^*} P_\chi$
for a monoid $P$ on $X/H$ with trivial stabilizer and $\nu_P=\nu_K$. Furthermore $K \cong \pi^*(P_\chi)$
for all $\chi \in H^*$. If the monoid $K$ has trivial 
stabilizer $H$, then  for any isogeny
$\pi:X\to Y$ the perverse sheaf $L=\pi_*(K)$ is an irreducible monoidal perverse sheaf on $Y$ with trivial stabilizer and the property $\nu_L=\nu_K$.
\end{Corollary}

{\it Proof}. Let $K$ be a monoid on $X$ with finite stabilizer $H$, $f:X\to Y$ be an isogeny with
$f(H)=0$.
Since $\pi$ is finite,  $L=\pi_*(K)=R\pi_*(K)$ is a semisimple perverse sheaf $L\neq 0$.
By corollary \ref{absch}, all summands $L_i$ of $L=\bigoplus_i P_i$ satisfy $\nu_{P_i} = \nu(K)$,
$\nu_L = \min(\nu_{P_i})= \nu_K$. By corollary \ref{mini}, 
at least one constituent $P (=Q_j)$ of $L$ is a monoid with $\nu_P=\nu_K$. 

\medskip
By the semisimplicity of $L$ and adjunction $$0\neq Hom_{\bf D}(P,L)\cong Hom_{\bf D}(P,\pi_*(K)) \cong Hom_{\bf D}(\pi^*(P),K) \ ,$$ Therefore there exists an exact sequence of perverse sheaves 
on $X$
$$0\to U\to \pi^*(P) \to K \to 0 \ $$
because any nontrivial morphism $\pi^*(P) \to K$ to the irreducible perverse sheaf $K$ is an epimorphism.
Since $\pi$ is finite, the functor $\pi_*$ is exact. Since $\pi_*\pi^*(P)\cong \bigoplus_{\chi\in Kern(\pi)^*} P_\chi$,
we get an exact sequence of perverse sheaves on $Y$  
$$0\to \pi_*(U)\to \bigoplus_{\chi\in Kern(\pi)^*} P_\chi \to L \to 0 \ .$$
Thus $L$ has at most $\# Kern(\pi)^*$ irreducible perverse constituents, and as twists of the monoid $P$ all of them are monoids of the same degree $\nu_K$. Hence the number of irreducible constituents of $L$ is $\dim({\cal H}^{-\nu_K}(L)_0)$. 
Since ${\cal H}^{-\nu_K}(L)_0 \cong \bigoplus_{x\in Kern(\pi)} {\cal H}^{-\nu_K}(K)_x $, 
in the case of the second assertion we get $\dim({\cal H}^{-\nu_K}(L)_0) \cong \dim({\cal H}^{-\nu_K}(K)_0) =1$ by
${\cal H}^{-\nu_K}(K)_x = 0$ for $x\neq 0$. So the second assertion follows immediately, since ${\cal H}^{-\nu_K}(L) \cong \delta_0$. For the first assertion the assumption $\# H = \# Kern(\pi)^*$ implies $\dim({\cal H}^{-\nu_K}(L)_0)=\# H$. Therefore $L$ has
$\#H$ irreducible constituents. Therefore $\pi_*(U)=0$, and  hence
$U=0$ and $K\cong \pi^*(P)$.  
\qed

\begin{Corollary}
Suppose $\pi: X\to Y$ is a separable isogeny and $K$ is an irreducible monoidal perverse sheaf on $Y$ with pullback $L=\pi^*(K)$. Then there exists an irreducible monoidal perverse sheaf $F$ with $\nu_F=\nu_K$ such that $L$ is isomorphic to the direct sum of translates $T^*_x(F)$, where $x$
runs over the cosets of $Kern(\pi)/Kern_F(\pi)$ for $Kern_F(\pi)=\{ x\in Kern(\pi)\ \vert
\ T_x^*(F)\cong F\} $. Furthermore $K_\chi \cong K$ holds for all $\chi$ whose pullback
$\chi \circ \pi_1(\tilde\pi)$ with respect to the isogeny $\tilde \pi: X/Kern_F(\pi) \to Y$ becomes trivial. 
\end{Corollary}

{\it Proof}. By etale descent one can show for an irreducible
perverse sheaf $K\in {\bf P}$ that the pullback $L=\pi^*(K)$ is a semisimple perverse sheaf and that the translations $T_x^*$ for $x\in Kern(\pi)$ act transitively on its simple constituents. Hence  $L=\bigoplus_i F_i$ for irreducible perverse sheaves $F_i$. Obviously $
\mu(K)=\nu_K\leq \mu(F_i)$. Notice that $\pi_*(L)=\pi_*\pi^*(K)= \bigoplus_\chi K_\chi$ implies $\pi_*(F_i) \cong K_\chi$ for some $\chi$. Since $\mu(F_i)=\mu(\pi_*(F_i))$, therefore $\nu_K\leq \mu(F_i) = \mu(\pi_*(F_i)) =
\mu(K_\chi)=\nu_K$ and this implies $\mu(F_i)=\nu_K$. Hence $\nu_{F_i} \leq \nu_K$.
But $\pi_*(F_i)\cong K_\chi$ 
implies $\nu_{F_i} = \nu_K$, by corollary \ref{absch}. 
Hence $\mu(F_i)=\nu_{F_i}=\nu_K$ for all
$i$. This shows that all $F_i$ are extremal and therefore $F_i \cong T^*_{x_i}(F)$ holds for certain $x_i\in X$, where $F$ is the unique constituent of $L=\pi^*(K)$ with the property ${\cal H}^{-\nu_K}(F)_0 \cong {\cal H}^{-\nu_K}(L)_0 \cong {\cal H}^{-\nu_K}(K)_0 \cong \Lambda$. In particular $F$ is a monoidal perverse sheaf on $X$ and $L$ is a direct sum of translates of $F$. This proves the first assertions.

\medskip
Since $F$ is invariant under translation by $Kern_F(\pi)$, $F$ descends to a perverse sheaf on $X/Kern_F(\pi)$ in the sense that $F\cong p^*(\tilde F)$ holds for $p:X\to X/Kern_F(\pi)$ and
$\tilde F$ is a constituent of $\tilde \pi^*(K)$. Then $Kern_{\tilde F}(\tilde \pi)=0$. We may therefore replace $\pi$ by $\tilde \pi$, 
So for the remaining statement we can assume $Kern_F(\pi)=0$ without restriction of generality. Then $L=\bigoplus_{x\in Kern(\pi)} T^*_x    
(F)$ and hence $\pi_*(L) = \#Kern(\pi) \cdot \pi_*(F)$. On the other hand $\pi_*(L)= \bigoplus_{\chi \in Kern(\pi)^*} K_\chi$. Both together imply that $K_\chi\cong K$ holds for all characters $\chi$  for  which $\chi \circ \pi_1(\tilde\pi)$ becomes trivial.
\qed

\medskip
By the adjunction formula $End(L) \cong Hom(K,\bigoplus_\chi K_\chi)$ for $L = \pi^*(K) \cong \bigoplus_x T_x^*(F)$,
we also conclude that $\#\{\chi\ \vert \ K\cong K_\chi \} \cdot \# \{ x\in Kern(\pi) \ \vert \ T_x^*(F) \cong F\}
= \# Kern(\pi)$. Here $\chi$ runs over all characters of $\pi_1(Y,0)$, whose restriction to $\pi_1(X,0)$ becomes trivial.

\medskip
{\bf Quasi-idempotents}. We either work in ${\bf D}$, or in  
a  hereditary localization $\bf D_H$ of ${\bf D}$ for
some hereditary class ${\bf H}$ as in [KrW], of course
possibly  $\bf D_H=D$.  Then $Hom_{\bf D_H}(P,{\bf D_H}^{> 0})=0$ for the image of some $P\in {\bf P}$ in $\bf D_H$. For the notation and further details we refer to [KrW].

\medskip
{\it Assumptions}. For a fixed  integer $d$, let $H^\bullet$ always denote graded $\Lambda$-vector spaces with the property $H^i=0$ for $\vert i\vert >d$. 
Suppose $$ P(X) \subset {\bf P} $$ is a class of simple objects 
closed under Tannaka duality, 
such that in $\bf D$
$$  K,L\in P(X) \ \Longrightarrow\  K*L \ \cong \  \bigoplus_{i\in I}\ H^\bullet(K,L,P_i) \otimes_\Lambda P_i \ \ \ \oplus \ \ T\ $$ for complexes $T$ in $\bf N_H$ and $P_i$ in $P(X)$. Here we assume that 
$P_i \not\cong P_j$ holds for $i \neq j$. By our assumption, $H^i(K,L,P_i)=0$ for $\vert i\vert > d$.

\begin{Lemma}\label{t1} Assume $P\in P(X)$ and ${\cal H}^{-d}(P)_0\neq 0$. Suppose $L[-d]\hookrightarrow K*P$
for $K,L\in P(X)$ but $K,L\not\in\bf N$. Then $K\cong L$. 
\end{Lemma}

{\it Proof}. By assumption $Hom_{\bf D_H}(K*P,L[-d])\neq 0$, and by rigidity this implies
$Hom_{\bf D_H}( P,K^\vee*L[-d])\neq 0$. Now, since $K^\vee*L[-d] =   \bigoplus_{i\in I} H^\bullet(K^\vee,L,P_i)\otimes_\Lambda P_i[-d]$ is in $ (\bigoplus_{i\in I} P_i) \oplus {\bf D_H}^{> 0}$ for some $I\subset {P}(X)$ (with multiplicities) again by our assumptions, we obtain
$ Hom_{\bf D_H}({P}, \ \bigoplus_{i\in I} P_i\ )\ \neq 0 $. Hence $ Hom_{\bf D_H}({P},P_i)\neq 0$ for some $i\in I$, and also $Hom_{\bf P}(P,P_i)\neq 0$ by [KrW, lemma 25]  for the simple objects  $P$ and $P_i$ in ${\bf P}$. So, $P_i\cong P$ are isomorphic as perverse sheaves. By the hard Lefschetz theorem, this defines in $\bf D$ a retract  
$ P[-d] \cong P_i[-d] \hookrightarrow H^{\bullet}(K^\vee,L,P_i)\otimes_\Lambda P_i\hookrightarrow K^\vee * L $.
Since ${\cal H}^0(P[-d])^*_0 \neq 0$, we get
$Hom_{\bf D}(L,K) = Hom_{\bf D}(K^\vee*L,\delta_0)= {\cal H}^0(K^\vee * L)^*_0 \neq 0$ and this implies $K\cong L$. 
\qed

\medskip
{\it For the next lemma} \ref{t2},  for arbitrary $K,P,L\in P(X)$ {\it we assume   in addition}: 
$$H^{-d}(K,P,L)=0 \ \Longrightarrow H^\bullet(K,P,L)=0\ .$$
 
\begin{Lemma} \label{t2} For $P\in P(X)$ assume ${\cal H}^{-d}(P)_0\neq 0$. 
Then for $K\in P(X)$ 
$$  K * P =   H^\bullet(K,P,K) \otimes_\Lambda K \quad \quad (\mbox{ in }  {\bf D_H}) \ .$$
 For monoids  $P\in P(X)$ not in $\bf N$ with  $\nu_P=d$, we get $ \dim_\Lambda(H^{-d}(P,P,P) )= 1$  and
 $$  P * P = H^\bullet(P,P,P) \otimes_\Lambda P \ ;$$
furthermore for $K\in P(X)$ either $H^\bullet(K,P,K)=H^\bullet(P,P,P)$ or $H^\bullet(K,P,K)=0$.
In particular, $P'*P=0$  holds in $\bf D_H$ for all monoids 
$P'\not\cong P$  with the property $\nu_{P'}=d$ under the assumption $P'\in P(X)$, but $P'\notin {\bf N}$.
\end{Lemma}

\medskip
{\it Proof}. If $H^\bullet(K,P,L_i)\neq 0$, then by our assumptions $H^{-d}(K,P,L_i)\neq 0$.
Hence $K*P= \bigoplus_{i\in I} H^\bullet(K,P,L_i) \otimes_\Lambda L_i$ for certain $L_i\in P(X)$ with $L_i[-d] \hookrightarrow K*{P}$. Hence $L_i\cong K$, by the last lemma. Since
$\dim_\Lambda(H^{-d}(P,P,P))$ for monoids $P={\cal P}$ with $\nu_P=d$
counts the multiplicity of  ${\cal P}[d]$ as a summand ${\cal P}[d] \hookrightarrow {\cal P}^{*2}$, this multiplicity is one by lemma 1.6) and 
${\cal P}_{\cal P} = {\cal P}$ (lemma 4). Now  $(H^\bullet(K,P,K)\otimes_\Lambda K)*P = (K*P)*P\cong K*(P*P)\cong H^\bullet(P,P,P) \otimes_\Lambda K*P$
in $\bf D_H$. For $K*P\neq 0$ this implies 
$H^\bullet(K,P,K)\cong H^\bullet(P,P,P)$. For
$K=P'$ and $P*P'\neq 0$ in $\bf D_H$, we get
$P*P' \cong H^\bullet(P,P',P)\otimes_\Lambda P = H^\bullet(P',P',P')\otimes_\Lambda P'$. Indeed, $P'$ satisfies the same conditions as $P$, so the roles of $P$ and $P'$ can be interchanged. A comparison in degree $-d$ gives $P\cong P'$.   \qed

\bigskip
{\bf Remark}.  In the above setting, $P\cong {\cal P}_K$ implies $H^\bullet(K,P,K)\neq 0$.

\goodbreak
\medskip
{\bf Quasi-idempotent complexes}. For given $L=\bigoplus_{i=-r}^r L_i[-i]$ with semisimple perverse sheaves $L_i$ for $-r \leq i \leq r$ assume
\begin{enumerate}
\item $L \cong L^\vee$ and $L_{-r} \cong L_r\neq 0$.
\item $L*L \cong H^\bullet \cdot L$ for $pr_{d}: H^{d}\cong \Lambda[-d]$ and $H^i=0$ for $\vert i\vert >d$.
\item There exists a commutative diagram with morphisms in the derived category
$$ \xymatrix{ L^\vee * L\ \  \ar[ddrr]_{eval_L} \ar[r]^\sim & \ \ L*L \ar[r]^\sim & H^\bullet \otimes_\Lambda L\ar[d]_{pr_d}\cr
& & H^{d}\otimes_\Lambda L \cong L[-d] \ar[d]\cr
& &  \delta_0} $$ 
\end{enumerate}
Then $H$ is selfdual and $H^{\pm d} \cong \Lambda$. Since ${\cal H}^\bullet(L^\vee*L)_0 \cong H^\bullet(X,D(L)\otimes_\Lambda^L L)$,
by condition 2 $$H^\bullet \otimes_\Lambda {\cal H}^\bullet(L)_0 \cong H^\bullet(X,D(L)\otimes_\Lambda L)\ $$ and both sides are independent of character twists, i.e. do not change when $L$ is replaced by $L_\chi$. Furthermore $L_r^\vee \cong L_{-r} \cong L_r$ by condition 1. Furthermore $L_r * L_r \cong L_r^\vee * L_r \neq 0$, since otherwise
the evaluation $eval _{L_r}=0$, and this implies $id_{L_r}=0$ and hence $L_r=0$ by rigidity. 

\medskip
Since $L_r*L_r \neq 0$, by the hard Lefschetz theorem ${}^pH^i(L*L)\neq 0 $ for some $i\geq 2r$. Hence $L*L \cong H^\bullet \otimes_\Lambda L$ implies  $2r \leq i \leq d+r$ or $r\leq d$.  
Let $\nu_L$ denote the minimum of all $\nu_C$  
for an irreducible perverse constituent
$C$ of some $L_i$. For the perverse amplitude $a(L_i,L_i)$ of $L_i*L_i$
and for $C\hookrightarrow L_i \cong L_i^\vee$ we have
$  \nu_C \leq a(L_i,L_i)$. Furthermore $a(L_i,L_i) + 2i \leq d+r$ by condition 2. Hence 
$\nu_C \leq a(L_i,L_i)\leq d+r -2i$. For $i=r$ this implies 
$\nu_C \ \leq\  d - r$, 
and hence  $\nu_L \leq \nu_C \leq d - r$. Therefore
$$   r  \ \leq \ d- \nu_L \ .$$
The morphism $eval_C=eval_{C[i]}$ for the direct summand $C^\vee * C = C[i]^\vee * C[i] \hookrightarrow L^\vee * L$
is obtained by restriction of  $eval_L$. Hence, by condition 3 there 
exists a commutative diagram
$$  \xymatrix{ {\cal P}_C[-\nu_C]\ar[rd]_{eval_C}  \ar[rr]^a &  & P_j[\lambda_j - d] \ar[ld]^b \cr
& \delta_0 & } $$
for some shifted irreducible summand $P_j \hookrightarrow L_{-\lambda_j}$ $$ P_j[\lambda_j] \hookrightarrow L \ .$$
Then $-\nu_C \leq \lambda_j - d$ and $\lambda_j - d \leq - \mu(P_j) \leq - \nu_{P_j}$, or otherwise
$a$ or $b$ is zero and hence $C=0$. For $b\neq 0$, also
$ \nu_{P_j} \leq \mu(P_j) \leq d - \lambda_j \leq \nu_C$.
For minimal $C$, i.e.  $\nu_C = \nu_L$, this implies
the equalities $\nu_{P_j} = \mu(P_j) = d - \lambda_j =\nu_C = \nu_L$.
The first equality gives ${\cal P}_C \cong P_j$, hence  $P_j$ is a minimal monoid. The last equality gives
$\lambda_j = d - \nu_L$, hence $r \leq \lambda_j$ 
from the inequality $r \leq d - \nu_L$ above. Since by our assumptions $\vert \lambda_j\vert \leq r$, therefore  $\lambda_j = r$ so that $P_j[\lambda_j] \hookrightarrow L_{-r}[r]$; in other words $$P_j \hookrightarrow L_{-r}\quad , \quad r + \nu_L = d \ . $$ 
Thus we found a multi-map from minimal constituents $C$ in $L$
to perverse minimal monoidal constituents $P_j$ in $L_{-r}$. On the other hand
$L_r^\vee \cong L_{-r} \cong L_r$, so for $C$ in $L_r\cong L_{-r} $ we get
${\cal P}_{C}[-\nu_C]\hookrightarrow C^\vee*C \hookrightarrow L_r*L_r$. 
On the other hand we found $\nu_C \leq a(L_r,L_r) \leq d-r$, so by the result
$d-r=\nu_L$ from above this implies $\nu_C\leq \nu_L$. Hence
all perverse constituents $C$ of $L_r$ are minimal, and  $$a(L_r,L_r)= \nu_L\ .$$
We claim that this implies that all perverse constituents of $L_{-r}$ (and hence 
of $L_r$) are monoids and that $L_{\pm r}$ is multiplicity free.

\begin{Proposition} \label{prop2}
Under the assumptions 1)-3) on $L= \bigoplus_{i=-r}^r L_i[-i]$, the two top and bottom perverse sheaves 
$L_r\cong L_{-r} = \bigoplus_j m_j \cdot K_j$ are multiplicity free perverse sheaves, i.e.
$m_j=1$ holds. Furthermore all the constituents $K_j$ are monoidal perverse sheaves
with $\nu_{K_j} = d - r = \nu_L$. 
\end{Proposition}

{\it Proof}. Recall $L_{-r}^\vee \cong L_r \cong L_{-r}$.
Therefore $m(K)=m({K^\vee})$ holds for the multiplicities $m(K)$ and $m(K^\vee)$
of $K$ and $K^\vee$ in $L_{-r}$. So,  for $m_j=m(K_j)$
$$ \bigoplus_j\ m_j^2 \cdot {\cal P}_{K_j}[\nu_{K_j}+2r]\hookrightarrow \bigoplus_j\ m_j^2 \cdot K_j^\vee[r]*K_j[r] \hookrightarrow
L_{-r}[r]*L_{-r}[r] \hookrightarrow L*L \cong H^\bullet\otimes_\Lambda  L \ .$$    
All $K_j$ in $L_{-r}$ are minimal, as shown already. Hence $\nu_{K_j} =\nu_L$, 
and $\nu_L+2r = d+r$ implies  
$$ \bigoplus_j\ m_j^2 \cdot {\cal P}_{K_j} \hookrightarrow H^\bullet[-d]\otimes_\Lambda L[-r]
 \ .$$    
Therefore $\bigoplus_j \ m_j^2 \cdot {\cal P}_{K_j} \hookrightarrow H^{-d}\otimes_\Lambda L_{-r} \cong L_{-r} $
and $\sum_{j, {\cal P}_{K_j}={\cal P}} m_j^2 \leq m({\cal P})$.
Since we already know that $K \in L_{\pm r}$ implies ${\cal P} = {\cal P }_K \in L_{\pm r}$, 
therefore $m({\cal P}) = 1$ follows and $m_j=0$  for all $K_j$ which are not monoids.
\qed

\medskip
{\bf Cohomology}. For irreducible $K\in Perv(X)$, define ${\cal S}(K)$ as the set of characters $\chi$
such that $H^\bullet(X,K_\chi)\neq H^0(X,K_\chi)$. For $\chi\in {\cal S}(K)$ define $h_\chi(K)$ to be the 
maximal
$i $ such that $H^{i}(X,K_\chi)\neq 0$. By the hard Lefschetz theorem 
$H^i(X,K_\chi)=0$ holds for $\vert i\vert > h_\chi(K)$ and $h_\chi(K)=h_\chi(K^\vee)\geq 0$.

\medskip
For $K\in E(X)$ the property  $H^\bullet(X,K_\chi)= H^0(X,K_\chi)$ is equivalent to $H^\bullet(X,K_\chi)=0$, using the preservation of the Euler  characteristic under character twists. Hence for $\nu_K>0$ this shows $\chi \in {\cal S}(K)$ iff $H^\bullet(X,K_\chi)\neq 0$. 
Therefore $ {\cal P}_K[\pm\nu_K] \hookrightarrow K*K^\vee$ and $K[\pm \nu_K]  \hookrightarrow K*{\cal P}_K$ imply 
$$ \fbox{$ \nu_K>0 \ \Longrightarrow \ {\cal S}({\cal P}_K)  =  {\cal S}(K) $} \ .$$
Furthermore, $h_\chi({\cal P}_K) +\nu_K \leq  h_\chi(K) + h_\chi(K^\vee)$ and 
$h_\chi(K) +\nu_K \leq h_\chi(K) + h_\chi({\cal P}_K)$ imply $$ \fbox{$ \nu_K\ \leq\ h_\chi({\cal P}_K )  \quad  \mbox{ for }\chi\in {\cal S}(K) $} \ .$$ Put $h_\chi(K) = \nu_K + e_\chi(K)$, then for all $ \chi\in {\cal S}(K)$ we obtain the inequalities 
$$ \fbox{$ 0 \ \leq \ e_\chi({\cal P}_K) \leq 2\cdot e_\chi(K)  $} \ .$$

\medskip
{\bf Relative case}. For a homomorphism $f:X\to Y$ we define $h_\chi^f(K)$, for all $\chi$ such that $Rf_*(K_\chi)\neq 0$,
to be the maximal integer $i$ for which ${}^p H^i(Rf_*(K_\chi))\neq 0$.
Since
$  H^k(X,K_\chi) \ =\ \bigoplus_{i+j=k} H^i(X,{}^p H^j(B,Rf_*(K_\chi))$
by the decomposition theorem, we obtain
$$   h_\chi(K) \ = \ \max_j \bigl(j + h({}^p H^{j}(Rf_*(K_\chi))\bigr) \ $$
where the maximum is taken now over all $j$ such that 
${}^p H^{j}(Rf_*(K_\chi))\neq 0$. Here we write  $h(F):=h_1(F)$ for the trivial character $\chi=1$. 
If $Rf_*(K_\chi)\neq 0$ is perverse, then $h_\chi(K)=h(Rf_*(K_\chi))$.
For all $\chi\in {\cal S}(K)$ 
$$ h_\chi(K)  \ \leq \ \dim(Kern(f)) + h_\chi\bigl( \bigoplus_i {}^pH^i(Rf_*(K_\chi))\bigr) \ .$$
%So let $\nu(X)$ be the maximum of $\nu(K), K\in F(X)$.
%We have $\nu(X) < \dim(X)$. And $K\boxtimes K$ and convolution gives
%for generic twists ${}^pH^i(Rf_*(K_\chi)) =0$ the maximum $i \leq \dim(A)$. 
%$$ 2 \nu_K \leq 2 h_\chi(K) \leq \dim(A)  + \nu_K \ .$$
%

\medskip
Let $P=K$ be an irreducible monoidal perverse sheaf on $X$ and $f:X\to Y$ be a homomorphism.
Then for every irreducible constituent $Q_j[\lambda_j]$ of $L=Rf_*(P)$ with perverse $Q_j$ we have
$h(P) \geq h(Q_j) + a_f(Q_j) \geq h(Q_j) + \lambda_j$, where $a_f(Q)= \max\{\lambda \vert  
Q[\lambda] \hookrightarrow L\}$ for a perverse sheaf $Q$.
On the other hand by theorem \ref{ungl} for every constituent $P_i[\lambda_i] \hookrightarrow L$
there exists some irreducible perverse sheaf $Q_j$ with $Q_j[\lambda_j] \hookrightarrow L$
and $\nu_P - \lambda_j \leq \nu_{P_i}$. For this particular $Q_j[\lambda_j]$
we conclude $\nu_P \leq \lambda_j + \nu_{P_i} $.
Now $-\nu_P \geq -\lambda_j - \nu_{P_i}$ together with  
$h(P) \geq h(Q_j) + \lambda_j$ gives the estimate $e(P)=h(P)-\nu_P \geq h(Q_j) - \nu_{P_i}$.
If $P_i$ is chosen minimal, then $Q_j \cong  {\cal P}_{P_i}$ by corollary \ref{mini}.
Therefore $e(Q_j)= h(Q_j) - \nu_{Q_j} = h(Q_j) - \nu_{P_i}$. So corollary \ref{mini} implies

\begin{Lemma}
For a monoidal perverse sheaf $P$ on $X$ and a homomorphism $f:X\to Y$
there exists a  monoidal perverse sheaf $Q$ on $Y$ such that
$Q[\nu_P - \nu_L] \hookrightarrow L=Rf_*(P)$ holds and  $e(Q) \leq e(P)$.
Furthermore $Q[\nu_P - \nu_L]$ is a minimal constituent of $Rf_*(P)$.
\end{Lemma}

\goodbreak

\bigskip\noindent
{\bf References}

\bigskip\noindent
[BBD] Belinson A., Bernstein J., Deligne P., {\it Faiscaux pervers}, Asterisque 100 (1982) 

\medskip\noindent
[D] Deligne P., {\it Categories Tensorielles}, Moscow Math. Journal, vol. 2, n.2, (2002), 227 - 248

%\medskip\noindent
%[FK] Franecki J., Kapranov M., {\it The Gauss map and a noncompact Riemann-Roch formula for constructible sheaves on semiabelian varieties}, arXiv:math/9909088

%\medskip\noindent
%[G] Ginsburg V., {\it Characteristic Varieties and vanishing cycles}, Invent. math. 84, 327 - 402 (1986) 

%\medskip\noindent
%[H] Hotta R., Takeuchi K., Tanisaki T., {\it D-Modules, Perverse sheaves, and representation theory}, Birkh\"auser Verlag (2008)
%
%\medskip\noindent
%[KFS] Kashiwara M., Fernandes T.M., Schapira P., {\it Truncated microspport and holomorphic solutions of $D$-modules}, arXiv:math/0203091
%
%\bigskip\noindent
%[KS] Kashiwara M., Shapira P., {\it Sheaves on manifolds}, Grundlehren der mathematischen Wissenschaften 292, Springer 2002
%
%
%\medskip\noindent
%[KS2] Kashiwara M., Shapira P., {\it Microlocal study of sheaves}, Asterisque, vol. 128, SMF Paris (1985)
%
%\medskip\noindent
%[K] Kashiwara M., {\it Systems of Microdifferential equations}, Progress in Math. 34. Birkh\"auser Verlag (1983)

\bigskip\noindent
[KrW] Kr\"amer T., Weissauer R., {\it Vanishing theorems for constructible sheaves on abelian varieties}, arXiv:1111.4947v2 (2011)

\bigskip\noindent
[BN] Weissauer R., {\it Brill-Noether sheaves}, arXiv:math/0610923v4 (2007)

\bigskip\noindent
[W2] Weissauer R., {\it A remark on the rigidity of BN-sheaves}, arXiv:1111.6095 
(2011)

\end{document}